\documentclass[a4paper, 12pt]{article}
\usepackage{amsmath}
\usepackage[all]{xy}
\usepackage{amsfonts} \usepackage{amssymb}\usepackage{latexsym}
\usepackage[english]{babel}
\usepackage{amscd}
\usepackage[dvips,final]{graphics}
\usepackage{locengli,math}
\usepackage{graphicx, Addaletree}

\begin{document}
\begin{center}
\textbf{\LARGE{\textsf{A rigidity theorem for $Moor$-bialgebras}}}
\vskip2cm
\large{Philippe {\sc Leroux}}\footnote{email: ph$\_$ler$\_$math@yahoo.com}
\end{center}
\vskip2cm

{\bf Abstract:}
We introduce the operad $Moor$, dual of the operad $NAP$ and the notion of $Moor$-bialgebras. We warn the reader that the compatibility relation linking the $Moor$-operation with the $Moor$-cooperation is not distributive in the sense of Loday.
Nevertheless,
a rigidity theorem (\`a la Hopf-Borel) for the category of connected $Moor$-bialgebras is given.
We show also that free permutative algebras can be equipped with a $Moor$-cooperation whose compatibility with the permutative product looks like the infinitesimal relation.

\noindent
\textbf{Notation}:
In the sequel $K$ is a characteristic zero field and $\Sigma_n$ is the group of permutations over $n$ elements. If $\mathcal{A}$ is an operad, then the $K$-vector space of $n$-ary operations is denoted as usual by $\mathcal{A}(n)$. We adopt Sweedler notation for the binary cooperation $\Delta$ on a $K$-vector space $V$ and set $\Delta(x)= x_{(1)} \otimes x_{(2)}$.

\section{Introduction}
The well-known Hopf-Borel theorem states that any connected cocommutative commutative bialgebra
(Hopf algebra) $\mathcal{H}$ is free and cofree over its primitive part $Prim \ \mathcal{H}$. Otherwise stated;
\begin{theo}(Hopf-Borel)
For any cocommutative commutative bialgebra $\mathcal{H}$ the following is equivalent.
\begin{enumerate}
 \item{$\mathcal{H}$ is connected;}
\item{$\mathcal{H}$ is isomorphic to $Com(Prim \ \mathcal{H})$ as a bialgebra;}
\item{$\mathcal{H}$ is isomorphic to $Com^c(Prim \ \mathcal{H})$ as a coalgebra.}
\end{enumerate}
\end{theo}
\noindent
In the theory developed by J.-L. Loday, this result is rephrased by saying that the triple of operads $(Com, Com, Vect)$ is good.
Some good triples of operads of type $(\mathcal{A},\mathcal{A}, Vect)$ or  $(\mathcal{C},\mathcal{A}, Vect)$ have been found since and a summary can be found in \cite{GB}.
The aim of this paper is to produce another good triple of operads of this form but without using the powerful theorems of J.-L. Loday \cite{GB}, simply because his first Hypothesis $(H0)$ is not fulfilled by our objects.

From our coalgebra framework on weighted directed graphs \cite{Coa, Ltrip}, we describe a directed graph by two cooperations $\Delta_M$ and $\tilde{\Delta}_M$ verifying:
$$(\tilde{\Delta}_M \otimes id)\Delta_M =(id \otimes \Delta_M)\tilde{\Delta}_M.$$
To code a bidirected graph, we have to add
the extra condition $\tau \Delta_M= \tilde{\Delta}$, where $\tau$ is the usual flip map. The previous equation becomes,
$$ (id \otimes \tau)(\Delta_M \otimes id)\Delta_M=(\Delta_M \otimes id)\Delta_M.$$
Such coalgebras were called $L$-cocommutative in \cite{Coa}. On the algebra side, this yield $K$-vector spaces equipped with a binary operation $\prec$ verifying,
$$ (x \prec y) \prec z = (x \prec z) \prec y.$$
Such algebras came out in the work of M. Livernet \cite{Liv} under the name nonassociative permutative algebras, $NAP$-algebras for short. The operad $NAP$ of $NAP$-algebras is important because it is related to the operad $preLie$ of preLie-algebras. Indeed, the triple
of operads $(NAP, preLie, Vect)$ has been shown to be good by M. Livernet \cite{Liv}. Requiring the operation $\prec$ to be associative leads to permutative algebras, or $Perm$-algebras for short \cite{Chp1}. In this paper, we introduce the dual, in the sense of Ginzburg and Kapranov \cite{GB}, of $NAP$-algebras, called $Moor$-algebras in Sections 2-3 and give a rigidity theorem for the category of connected $Moor$-bialgebras in Section 4, that is the triple of operads $(Moor, Moor, Vect)$ is good. This category is interesting, as we said, for we cannot apply the powerful results of J.-L. Loday \cite{GB} since the compatibility relation linking the cooperation and the operation of a $Moor$-bialgebra is not distributive as required in \cite{GB}, Hypothesis $(H0)$. We end with Section 5, where we show that the free permutative algebra over a $K$-vector space $V$ can be equipped with a $Moor$-cooperation whose compatibility relation with the permutative product looks like the nonunital infinitesimal relation.

\section{On $Moor$-algebras}
Define the operad $Moor$, ($Moor$ because a typical element of a $Moor$-algebra looks like $(\ldots((x_1x_2)x_3) \ldots)x_n)$ whose parentheses are concentrating at the beginning, reminding boats being moored one behind the other) to be the free operad on one binary operation $\prec$
divided out by the following set of relations:
$$ R:= \{ (x \prec y) \prec z - (x \prec z) \prec y; \ \ x \prec (y \prec z) \}.$$
If $V$ stands for a $K$-vector space, then
$S(V)$ stands for the symmetric module over $V$, that is:
$$ S(V):= K \oplus \bigoplus_{n>0} \ S^n(V), $$
where $S^n(V)$ is the quotient of $V^{\otimes n}$ by the usual action of the symmetric group $\Sigma_n$. A typical element of $S^n(V)$ will be written $v_1 \vee v_2 \vee \ldots \vee v_n,$
where the $v_i \in V$.
\begin{theo} The following hold.
\begin{enumerate}
\item {The dual of the operad $NAP$ is the operad $Moor$.}
\item {The free $Moor$-algebra over a $K$-vector space $V$ is $V \otimes S(V)$ as a $K$-vector space equipped with the operation $\prec$ defined by:
$$ v \otimes \omega \prec v' \otimes \omega' = v \otimes \omega \vee v',$$
if $\omega' \in K$ and vanishes otherwise.}
\item {The generating series of the operad $Moor$ is,
$$ f_{_{Moor}}(x):=xe^{x}= \sum_{n>0} \  n \frac{x^n}{n!}.$$}
\end{enumerate}
\end{theo}
\Proof
Observe that the free binary operad $\mathcal{F}$ with one binary operation obey the relation $\dim \ \mathcal{F}(3)=12$. We get $\dim \ NAP(3)=9$ and $\dim \ Moor(3)=3.$ As in $\mathcal{F}(3)$, quadratic relations defining $NAP(3)$ are orthogonal (see \cite{GK}) to those defining $Moor(3)$, the dual of $NAP$ is $Moor$. Let $V$ be a $K$-vector space.
The $K$-vector space $V \otimes S(V)$  equipped with the operation,
$$ \prec: \ V \otimes S(V) \bigotimes V \otimes S(V) \rightarrow V \otimes S(V), \ \ \ v \otimes \omega \prec v' \otimes \omega' = v \otimes \omega\vee v',$$
if $\omega' \in K$ and vanishes otherwise is a $Moor$-algebra.
Observe that $i: V \hookrightarrow V \otimes K \hookrightarrow V \otimes S(V)$ defined by $i(v):= v \otimes 1_K$ realises the expected embedding.
Let $(A, \prec_A)$ be a $Moor$-algebra and $f: V  \rightarrow A$ be a map. Define $\tilde{f}: V\otimes S(V) \rightarrow A$ by,
$$ \tilde{f}(v \otimes 1_K):=f(v),$$
$$ \tilde{f}(v \otimes v_1 \cdots v_p):=(\cdots ((f(v)\prec_A f(v_1))\prec_A f(v_2)) \ldots \prec_A f(v_{p-1})) \prec_A f(v_p)).$$
Then, $\tilde{f}$ is a $Moor$-morphism and the only one such that $\tilde{f} \circ i =f.$
For the last item,
observe that in a $Moor$-algebra only these monomials,
$$ (left \ combs: \ (lc)) \ \ \ (\cdots (v_1\prec v_2)\prec v_3) \ldots \prec v_{n-1}) \prec v_n),$$
do not vanish. Indeed, one can model $n$-ary operations of the $Moor$-operad with planar rooted binary trees whose nodes are decorated by $\prec$. For instance, $x \prec(y \prec z)$ is represented by $\treeBA_{\prec}$ so $\treeBA_{\prec}=0$ and only left combs survive. Therefore, we get $n(n-1)!$ such left combs but because of the
relation $(x \prec y) \prec z = (x \prec z) \prec y$, the relation $(lc)$ is invariant under the action of the symmetric group $\Sigma_{n-1}$.
Hence, $\dim Moor(n)=n$.
\eproof

\section{The cofree $Moor$-coalgebra}
\noindent
Let $i$ be an integer. By $v^{i}$, we mean $v \vee \ldots \vee v$, times $i$. In the sequel, we set by induction, for all $n>0$,
$\Delta_{_{\mathcal{H}}}^{(1)}=\Delta_{_{\mathcal{H}}}$ and
$\Delta_{_{\mathcal{H}}}^{(n)}:=(\Delta_{_{\mathcal{H}}} \otimes id_{ (n-1)} )\Delta_{_{\mathcal{H}}}^{(n-1)}$ for any cooperation
$\Delta_{_{\mathcal{H}}}$ of a $Moor$-coalgebra $\mathcal{H}$.
We get the following two propositions by dualising the corresponding results in the proof of Theorem~2.1.

\begin{lemm}
\label{sigma}
Let $(\mathcal{H}, \Delta_{_\mathcal{H}})$ be a coalgebra whose cooperation verifies
$\Delta_{_\mathcal{H}}^{(2)}=(id \otimes \tau)\Delta_{_\mathcal{H}}^{(2)}$.
For all $n>0$,  the map: $$\Delta_{_\mathcal{H}}^{(n)}: \mathcal{H} \rightarrow \mathcal{H}^{\otimes (n+1)}, \ x \mapsto \Delta^{(n)}(x):=x_{n+1} \otimes x_{n} \otimes \ldots \otimes x_i \otimes \ldots \otimes x_2 \otimes x_1,$$
has its last $n$ components invariant by $\Sigma_{n}$.
\end{lemm}
\Proof
Fix $i=1, \ldots, n-1$. The following,
\begin{eqnarray*}
\Delta_{_\mathcal{H}}^{(n)} &=& (\Delta_{_\mathcal{H}}^{(n-i-1)}\otimes id_{ (i+1)})\circ (\Delta_{_\mathcal{H}}^{(2)} \otimes id_{ (i-1)})\circ \Delta_{_\mathcal{H}}^{(i-1)},\\
&=& (\Delta_{_\mathcal{H}}^{(n-i-1)}\otimes id_{ (i+1)})\circ ((id \otimes \tau)\Delta_{_\mathcal{H}}^{(2)} \otimes id_{ (i-1)})\circ \Delta_{_\mathcal{H}}^{(i-1)},\\
&=& (id_{(n-i)}\otimes \tau \otimes id_{(i-1)})\circ (\Delta_{_\mathcal{H}}^{(n-i-1)}\otimes id_{ (i+1)})\circ (\Delta_{_\mathcal{H}}^{(2)} \otimes id_{ (i-1)})\circ \Delta_{_\mathcal{H}}^{(i-1)},
\end{eqnarray*}
shows that the last $n$ components of $\Delta_{_\mathcal{H}}^{(n)}$ are invariant by the transpositions $(i, i+1)$ for all $i=1, \ldots, n-1,$ hence the claim.
\eproof

\begin{prop}
The cofree $Moor$-coalgebra over a $K$-vector space $V$  is:
$$ Moor^c(V):= V \otimes S(V),$$
as a $K$-vector space equipped with the following
co-operation $\delta$ defined as follows:
$$ \delta(v\otimes 1_K)=0,$$
$$ \delta(v_1\otimes v_2^{i_2} \vee \ldots \vee v_n^{i_n})=\sum_{k=2}^n (v_1 \otimes v_{2}^{i_2}\vee \ldots \vee v_{k}^{i_k-1}\vee \ldots \vee v_{n}^{i_n})\otimes (v_{k}\otimes 1_K).$$
\end{prop}

\noindent
Let $\Gamma  V^{\otimes n}$ be the $K$-vector space of tensors invariant through the usual action of $\Sigma_n$. For all $n>1$, define $j_n: V\otimes \Gamma  V^{\otimes n} \rightarrow V\otimes S^n(V)$ by
$j_n(\sum_{\sigma \in \Sigma_n} \ v \otimes v_{1} \otimes \ldots \otimes v_{n})=v \otimes v_{1} \vee \ldots \vee v_{n}$. They are bijective maps since $K$ is a characteristic zero field.
\begin{prop}
\label{exten}
If $(\mathcal{H},\Delta_{_\mathcal{H}})$ is a $Moor^c$-coalgebra and $f:\mathcal{H} \rightarrow V$ a linear map, set by induction $f^{\otimes 1}=f$ and $f^{\otimes n}=f^{\otimes (n-1)} \otimes f $.
Then, the map $\tilde{f}:\mathcal{H} \rightarrow Moor^c(V)$
defined by:
$$ \tilde{f}:= \sum_{n=1}^\infty \ j_n \circ f^{\otimes (n+1)} \circ \Delta_{_\mathcal{H}}^{(n)},$$
is the unique coalgebra morphism verifying $\pi \circ \tilde{f}=f$, where $\pi: Moor^c(V) \twoheadrightarrow V$ is the canonical projection.
\end{prop}

\section{On $Moor$-bialgebras}
\subsection{Definition}
\noindent
In the sequel, we set for any $v_1, \ldots, v_n \in V$:
$$(\ldots(v_1 \prec v_2)\prec \ldots) \prec v_n):=[v_1|v_2, \ldots, v_n].$$

\noindent
By definition, a $Moor$-bialgebra $\mathcal{H}$ is the data of:
\begin{enumerate}
\item {A graduated $Moor$-algebra  $\mathcal{H}:=\bigoplus_{p>0} \ \mathcal{H}_p$,}
\item{A $Moor$-cooperation $\Delta_{_\mathcal{H}}: \mathcal{H}\rightarrow \mathcal{H}^{\otimes 2}$, i.e., verifying:
$$ (id  \otimes \Delta_{_\mathcal{H}})\Delta_{_\mathcal{H}}=0,$$
$$ (\Delta_{_\mathcal{H}} \otimes id)\Delta_{_\mathcal{H}}= (id \otimes \tau)(\Delta_{_\mathcal{H}} \otimes id)\Delta_{_\mathcal{H}},$$}
\item{The $Moor$-operation and cooperation have to be related by the following compatibility condition:
$$ \Delta_{_\mathcal{H}}(x \prec y)= x \otimes e(y) + (x_{(1)}\prec y)\otimes x_{(2)},$$
for any $x,y \in \mathcal{H}$, where $e: \mathcal{H}\twoheadrightarrow \mathcal{H}_1$ is the canonical projection.}
\end{enumerate}
A morphism of $Moor$-bialgebras is a morphism of graduated $Moor$-algebras and a morphism of $Moor$-coalgebras.
Observe that this compatibility relation is not distributive in the sense of J.-L. Loday \cite{GB}.
By $Prim \ \mathcal{H}:= \ker \Delta_{_\mathcal{H}}$ we mean the $K$-vector space of primitive elements.

\begin{prop}
Let $\mathcal{H}$ be a $Moor$-bialgebra. Then, $\ker \Delta_{_\mathcal{H}} = \tilde{\mathcal{H}_1}  \oplus \tilde{\mathcal{H}},$
where $\tilde{\mathcal{H}}:=\bigoplus_{j \in J} \ \tilde{\mathcal{H}_j},$
with $J$ a suitable subset of $\mathbb{N}\setminus \{0,1\}$. If it exists, $\tilde{\mathcal{H}}$ is a $Moor$-algebra equipped with the following action:
$$ \tilde{\mathcal{H}_1} \otimes \tilde{\mathcal{H}} \rightarrow \tilde{\mathcal{H}}, \ \ h_1 \otimes h \mapsto h_1\prec h.$$
Moreover, $\mathcal{H}_*:=\bigoplus_{p>1} \ \mathcal{H}_p$ acts on $\ker \Delta_{_\mathcal{H}}$ on the right via:
$$ \ker \Delta_{_\mathcal{H}} \otimes  \mathcal{H}_* \rightarrow \ker \Delta_{_\mathcal{H}}, \ \ \ a \otimes h \mapsto a \prec h.$$
\end{prop}
\Proof
For $j>0$, set $\tilde{\mathcal{H}_j}$ the $K$-vector space of primitive elements of degree $j$. For $h_j$ and $h_{j'}$ two primitive elements of degrees resp. $j\geq 1$ and $j'>1$, one has:
$$\Delta_{_\mathcal{H}}(h_j \prec h_ {j'})=h_j \otimes e(h_ {j'})=0,$$
hence $h_j \prec h_ {j'}$ of degree $j+j'$ is primitive. If $h \in \mathcal{H}_*$, then
$\Delta(h_j \prec h)= h_j \otimes e(h)=0$ hence the last claim.
\eproof

\noindent
Set $\Delta^{(n)}_{_\mathcal{H}}:= (\Delta_{_\mathcal{H}} \otimes id_{(n-1)})\Delta^{(n-1)}_{_\mathcal{H}}$ with
$id_{(n-1)}= \underbrace{id \otimes \ldots \otimes id}_{times \ n-1}$ for all $n\geq 1$.
By definition, a $Moor$-bialgebra is is said to be connected if $\mathcal{H}=\cup_{r>0} \ \mathcal{F}_r\mathcal{H},$
where the filtration $(\mathcal{F}_r\mathcal{H})_{r>0}$ is defined as follows:
$$ (The \ primitive \ elements:) \ Prim \ \mathcal{H}:= \mathcal{F}_1\mathcal{H}:= \ker \Delta_{_\mathcal{H}} \subset \mathcal{H}_1.$$
Set $\Delta^{(n)}_{_\mathcal{H}}:= (\Delta_{_\mathcal{H}} \otimes id_{(n-1)})\Delta^{(n-1)}_{_\mathcal{H}}$ with
$id_{(n-1)}= \underbrace{id \otimes \ldots \otimes id}_{times \ n-1}$ for all $n\geq 1$. Then,
$$ F_r\mathcal{H}:=\ker \ \Delta^{(r)}_{_\mathcal{H}}.$$
\noindent
Here is an example of connected $Moor$-bialgebras.

\begin{theo}
Let $V$ be a $K$-vector space.
The free $Moor$-algebra over $V$ is a connected  $Moor$-bialgebra.
\end{theo}
\Proof
Let $V$ be a $K$-vector space.
Define the co-operation $\Delta$ by induction as follows:
$$ \Delta(v \otimes 1_K):=0,$$
$$ \Delta(x \prec y)= x \otimes \pi(y) + (x_{(1)}\prec y)\otimes x_{(2)},$$
for any $v \in V$, $x,y \in Moor(V)$, where $\pi: Moor(V)\twoheadrightarrow i(V)$ is the canonical projection map. As $x \prec y=0$ for all $x,y \in Moor(V)$ and $y$ of degree at least 2, we have to check that $\Delta (x \prec y)$ vanishes. But,
$$
\Delta (x \prec y) = x \otimes \pi(y) + (x_{(1)}\prec y)\otimes x_{(2)}=0,
$$
because $\pi(y)=0$ since the degree of $y$ is at least 2 and  $x_{(1)}\prec y=0$ for the same reason.
By induction one proves:
$$ \Delta([v_1|v_2, \ldots, v_n])=\sum_{k=2}^n \ [v_1|v_2, \ldots, \hat{v}_k, \ldots, v_n]\otimes (v_k\otimes 1_K),$$
for all $v_1, \ldots, v_n \in V$ and where the hat notation means as usual that the involved element vanishes. From that formula, it is straightforward to check that the co-operation $\Delta$ verifies:
$$ (id  \otimes \Delta)\Delta=0,$$
$$ (\Delta \otimes id)\Delta= (id \otimes \tau)(\Delta \otimes id)\Delta.$$
Hence, the free $Moor$-algebra over $V$, which is graduated by construction, is a $Moor$-bialgebra.

The map $\phi(V): Moor(V) \rightarrow Moor^c(V)$, defined as follows:
$$ \phi(V)(v \otimes 1_K)=v\otimes 1_K,$$
$$ \phi(V)([v_1|v_{2}^{i_2}, \ldots, v_{n}^{i_n}])=i_2!\ldots i_n! \ v_1\otimes v_{2}^{i_2}\vee \ldots \vee v_{n}^{i_n},$$
is an isomorphism of $Moor$-coalgebras.
It suffices to observe that:
$$\delta(v_1\otimes v_2^{i_2} \vee \ldots \vee v_n^{i_n})=\sum_{k=2}^n (v_1 \otimes v_{2}^{i_2}\vee \ldots \vee v_{k}^{i_k-1}\vee \ldots \vee v_{n}^{i_n})\otimes (v_{k}\otimes 1_K),$$
and:
$$\Delta([v_1| v_2^{i_2} , \ldots , v_n^{i_n}])=\sum_{k=2}^n i_k[v_1 | v_{2}^{i_2}\vee \ldots \vee v_{k}^{i_k-1}\vee \ldots \vee v_{n}^{i_n}]\otimes (v_{k}\otimes 1_K),$$
Thus:
$$\delta (\phi(V)([v_1| v_2^{i_2} , \ldots , v_n^{i_n}]))=i_2!\ldots i_n! \sum_{k=2}^n (v_1 \otimes  v_{2}^{i_2}\vee \ldots \vee v_{k}^{i_k-1}\vee \ldots \vee v_{n}^{i_n})\otimes (v_{k}\otimes 1_K),$$
and:
\begin{eqnarray*}
& & (\phi(V) \otimes \phi(V))\Delta ([v_1| v_2^{i_2} , \ldots , v_n^{i_n}]) =\\
&=& (\phi(V) \otimes \phi(V)) (\sum_{k=2}^n i_k \ [v_1 |  v_{2}^{i_2}, \ldots,  v_{k}^{i_k-1},\ldots, v_{n}^{i_n}]\otimes (v_{k}\otimes 1_K))\\
&=& \sum_{k=2}^n i_2!\ldots i_k(i_k-1)! \ldots i_n! \ (v_1 \otimes v_{2}^{i_2}\vee \ldots \vee v_{k}^{i_k-1}\vee \ldots \vee v_{n}^{i_n})\otimes (v_{k}\otimes 1_K).
\end{eqnarray*}
Hence, $\phi(V)$ is a coalgebra morphism and is bijective since $K$ is a characteristic zero field. Therefore, $\ker \Delta = (Moor(V))_1$ and the filtration being given by the $((Moor(V))_n)_{n>0}$, the free $Moor$-algebra over $V$ is a connected $Moor$-bialgebra.
\eproof

\begin{lemm}
\label{prim}
A connected $Moor$-bialgebra $\mathcal{H}$ is generated by its primitive elements. Moreover, $\ker \Delta_{_\mathcal{H}} = \mathcal{H}_1$.
\end{lemm}
\Proof
Let $x \in \mathcal{F}_r \mathcal{H}$
with $r$ minimal and belongs to $\mathcal{H}_p$, $p>0$ which is not primitive. Write $\Delta_{_\mathcal{H}}(x)=x_{(1)}\otimes x_{(2)}$ as a sum of independents vectors.
We get $0=(id \otimes \Delta_{_\mathcal{H}})\Delta_{_\mathcal{H}}(x)= x_{(1)}\otimes \Delta_{_\mathcal{H}}(x_{(2)})$. Hence $\Delta_{_\mathcal{H}}(x_{(2)})=0$ and the $x_{(2)}$
are primitive elements and belongs to $\mathcal{H}_1$. Moreover,
$0=\Delta_{_\mathcal{H}}^{(r)}(x)= \Delta_{_\mathcal{H}}^{(r-1)}(x_{(1)})\otimes x_{(2)}$ which leads to $\Delta_{_\mathcal{H}}^{(r-1)}(x_{(1)})=0$ and the $
x_{(1)} \in \mathcal{F}_{r-1} \mathcal{H}$. Therefore,
$$\Delta_{_\mathcal{H}}^{(r-1)}(x)=x_{(1)}\otimes \ldots \otimes x_{(r)},$$
where the $x_{(i)}$ for $1\leq i \leq r$ are primitive. However,
$$ \Delta_{_\mathcal{H}}^{(r-1)}(x-[x_{(1)}|x_{(2)}, \ldots, x_{(r)}])=0,$$
hence $x-[x_{(1)}|x_{(2)}, \ldots, x_{(r)}]$ is a primitive element and belongs to $\mathcal{H}_1$. As $x \in \mathcal{H}_p$, $[x_{(1)}|x_{(2)}, \ldots, x_{(r)}]\in \mathcal{H}_r$, we get
$p=r$ and $x=[x_{(1)}|x_{(2)}, \ldots, x_{(r)}]$.
\eproof

\subsection{A rigidity theorem for connected $Moor$-bialgebras }
\begin{theo}
A connected $Moor$-bialgebra $ \mathcal{H}$ is free and cofree over its primitive part $Prim \ \mathcal{H}$, that is the following is equivalent for any $Moor$-bialgebra $\mathcal{H}$:
\begin{enumerate}
 \item{$\mathcal{H}$ is connected;}
\item{$\mathcal{H}$ is isomorphic to $Moor(Prim \ \mathcal{H})$ as a $Moor$-bialgebra;}
\item{$\mathcal{H}$ is isomorphic to $Moor^c(Prim \ \mathcal{H})$ as a $Moor$-coalgebra.}
\end{enumerate}
\end{theo}
\Proof
Let $\mathcal{H}$ be a connected $Moor$-bialgebra. Since, $Moor(Prim \ \mathcal{H})$ is free, we get:
$$\xymatrix{
Prim \ \mathcal{H} \ar[r]^-{i} \ar@{->}[rd]_{j} & Moor(Prim \ \mathcal{H}) \ar@{->}[d]^{\tilde{i}} \\
 & \mathcal{H}}$$
where $\tilde{i}$ is the unique $Moor$-morphism verifying $\tilde{i}\circ i=j$, where $i$ and $j$ are the canonical injections. Via Lemma \ref{prim}, $\tilde{i}$ is surjective.
Since $Moor^c(Prim \ \mathcal{H})$ is cofree, we get,
$$\xymatrix{
\mathcal{H} \ar[r]^-{\tilde{e}} \ar@{->}[rd]_{e} & Moor^c(Prim \ \mathcal{H}) \ar@{->>}[d]^{\pi} \\
 & Prim \ \mathcal{H}}$$
with $\tilde{e}$ the unique morphism of coalgebra extending the canonical projection $e$.
Still set by induction
$\Delta_{_{\mathcal{H}}}^{(1)}=\Delta_{_{\mathcal{H}}}$ and
$\Delta_{_{\mathcal{H}}}^{(n)}:=(\Delta_{_{\mathcal{H}}} \otimes id_{ (n-1)} )\Delta_{_{\mathcal{H}}}^{(n-1)}$. Set $e^{\otimes 1}=e$ and $e^{\otimes n}=e^{\otimes (n-1)} \otimes e $.
Recall the coalgebraic morphism $\tilde{e}$ is given as follows:
$$ \tilde{e}(x)=\sum_{n=1}^\infty \ j_n \circ e^{\otimes (n+1)}\circ \Delta_{_{\mathcal{H}}}^{(n)}(x).$$
As a connected $Moor$-bialgebra is generated by its primitive elements, we focus on elements $x:=[x_1| x_2, \ldots, x_n]$, with the $x_i$ primitive. As expected,
$$ \tilde{e}([x_1| x_2^{i_2}, \ldots, x_n^{i_n}]_\prec)= i_2! \ldots i_n! \ x_1\otimes x_2^{i_2}\vee \ldots \vee x_n^{i_n}.$$
Hence, we get on the whole $\mathcal{H}$, $\phi(Prim \ \mathcal{H})=\tilde{e}\circ \tilde{i}$ where $\phi(Prim \ \mathcal{H})$ is defined in the proof of Theorem 4.2. Since $\tilde{i}$ is surjective (Lemma \ref{prim})
and $\phi(Prim \ \mathcal{H})$ is bijective, $\tilde{i}$ is injective and is an isomorphism. Hence $\tilde{e}$ is also an isomorphism of $Moor$-coalgebras since $\tilde{e}= \phi(Prim \ \mathcal{H})\circ\tilde{i}^{-1}$.
\eproof

\section{A $Moor$-cooperation over free $Perm$-algebras}
Permutative algebras have been introduced in \cite{Chp1}. In fact the following holds.
\begin{prop}
Let $V$ be a $K$-vector space. Then, the $K$-vector space $V\otimes S(V)$ equipped with the operation $\sqsupset$ defined by,
$$ v_1 \otimes v_2 \vee \ldots \vee v_n \sqsupset w_1\otimes w_2 \vee \ldots \vee w_m = v_1 \otimes v_2 \vee \ldots \vee v_n \vee w_1 \vee w_2 \vee \ldots \vee w_m,$$
for all $v_i,w_j \in V$,
is the free $Perm$-algebra over $V$.
\end{prop}

\noindent
A $Perm$-algebra $P$ is said to be unital if it exits an element denoted by 1 such that
$ x \sqsupset 1 = x,$
for all $x \in P$, the symbols $ 1 \sqsupset x$, $1 \sqsupset 1$ being not defined. The augmented free $Perm$-algebra over a $K$-vector space $V$, $K \oplus Perm(V)$, is a unital $Perm$-algebra.

\noindent
Let $V$ be a $K$-vector space. On $V \otimes S(V)$, one can
define the left and right maps as follows:
$$ l(1_K)=0, \ \ \ l(v_1\otimes v_2 \vee \ldots \vee v_n)=v_1\otimes 1_K, \ \ \ \ \ r(v_1 \otimes v_2\vee \ldots\vee v_n)= \frac{1}{n-1}\sum_{i=2}^{n} \ v_i \otimes v_2 \vee \ldots \vee \hat{v}_i\vee \ldots\vee v_n,$$
$$ r(v_1 \otimes 1_K)=1_K, \ \ \ r(1_K)=0,$$
for all $v_i \in V$.
\begin{prop}
Let $V$ be a $K$-vector space. Then, the augmented free $Perm$-algebra over a $K$-vector space $V$, $K \oplus Perm(V)$, can be equipped with a $Moor$-cooperation
$\Delta$ verifying the following compatibility relation:
$$ \Delta(1_K)=0.$$
$$ \Delta(x \sqsupset y)= (x_{(1)}\sqsupset y) \otimes x_{(2)} + (x\sqsupset y_{(1)}) \otimes y_{(2)}+ (x \sqsupset r(y)) \otimes l(y),$$
for all $x,y \in K\oplus Perm(V)$.
\end{prop}
\Proof
Recall in $V \otimes S(V)$ the existence of the following $Moor$-cooperation $\Delta$:
$$ \Delta(v_1 \otimes v_2 \vee \ldots \vee  v_n)=\sum_{k=2}^n \ (v_1 \otimes v_2\vee \ldots \vee \hat{v}_k\vee \ldots\vee v_n) \otimes (v_k\otimes 1_K),$$
defined for all $v_i \in V$.
Add $\Delta(1_K):=0$.
Hence, set $x=v_1 \otimes v_2 \vee \ldots \vee  v_n$ and $y=w_1 \otimes w_2 \vee \ldots \vee  w_m$ and observe that
$$ \Delta(x \sqsupset y)= (x_{(1)}\sqsupset y) \otimes x_{(2)} + (x\sqsupset y_{(1)}) \otimes y_{(2)}+ (x \sqsupset r(y)) \otimes l(y),$$
holds.
\eproof

\noindent
\textbf{Acknowledgments:}
Many thanks to M. Livernet and J.-L. Loday for usefull discussions.
\bibliographystyle{plain}
\bibliography{These}

\end{document}